\documentclass{article}
\usepackage{amssymb}
\usepackage{srcltx}
\usepackage{amsfonts}
\usepackage{amsmath}

\begin{document}

\begin{center}
{\Large Some identities in algebras obtained by the Cayley-Dickson process}

\begin{equation*}
\end{equation*}

Cristina FLAUT \ and $\ $Vitalii \ SHPAKIVSKYI \ {\large \ }%
\begin{equation*}
\end{equation*}
\end{center}

\textbf{Abstract. }{\small In this paper we will prove that the Hall
identity is true in all algebras obtained by the Cayley-Dickson process and,
in some conditions, the converse is also true. }%
\begin{equation*}
\end{equation*}


\textbf{2000} \textbf{AMS}{\small \ Classification: \ 15A66, 17A05, 17A20,
17A35, 17A45. \ }%
\begin{equation*}
\end{equation*}

\textbf{0.} \textbf{Introduction}

\begin{equation*}
\end{equation*}

In October 1843, William Rowan Hamilton discovered the quaternions, a $4$%
-dimensional algebra over $\mathbb{R}$ which is associative and
noncommutative algebra. In December $1843,$ \ John Graves discovered the
octonions, an $8$-dimensional algebra over $\mathbb{R}$ which is
nonassociative and noncommutative algebra. These algebras were rediscovered
\ by Arthur Cayley in $1845$ and are also known sometimes as the Cayley
numbers. This process, of passing from $\mathbb{R}$ to $\mathbb{C}$, from $%
\mathbb{C}$ to $\mathbb{H}$ and from$\ \mathbb{H}$ to $\mathbb{O}$ has been
generalized to algebras over fields and over rings. It is called the \ 
\textit{Cayley-Dickson doubling process} \ or the \textit{Cayley--Dickson
process. } In 1878, W. K. Clifford discovered Clifford algebras defined to
have generators $e_{1},$ $e_{2},...,e_{n}$ which anti-commute and satisfy $%
e_{i}^{2}=a_{i}\in \mathbb{R},$ for all $i\in \{1,2,...,n\}.$ These algebras
generalize the real numbers, complex numbers and quaternions( see [Le; 06 ])

Even if are old, quaternions, octonions and Clifford algebras have at
present many applications, as for example in physics, coding theory,
computer vision, etc. For this reasons these algebras are intense studied.
In [Ha; 43], Hall proved that the identity \ $\left( xy-yx\right)
^{2}z=z\left( xy-yx\right) ^{2}$\ \ holds for all elements $x,y,z$ in a
quaternion algebra.\ This identity is called \textit{Hall identity}. \
Moreover, he also proved the converse: if the Hall identity is true in a
skew-field $F,\ $then $F$ is a quaternion division algebra. In [Smi; 50],
Smiley proved that the Hall identity is true for the octonions and he also
proved the converse: if the Hall identity is true in an alternative division
algebra $A,\ $then $A$ \ is an \ octonion division algebra.

In this paper we will prove that the Hall identity is true in all algebras
obtained by the Cayley-Dickson process and, in some conditions, the converse
is true for split quaternion algebras and split octonion algebras. 
\begin{equation*}
\end{equation*}

\textbf{1. Preliminaries}%
\begin{equation*}
\end{equation*}

In this paper, we assume that $K$ is a commutative field with $charK\neq 2$
and $A$ is an algebra over the field \ $K.$ The \textit{center} \ $C$ of an
algebra $A$ is the set of all elements $c\in A$ which commute and associate
with all elements $x\in $ $A.$ An algebra $A$ is a \textit{simple} algebra
if \ $A$ is not a zero algebra and \ $\{0\}$ and $A$ are the only ideals of $%
A.$ The algebra $A$ is called \textit{central simple} if \ the algebra\ $%
A_{F}=F\otimes _{K}A$ is simple for every extension $F$ of $K.$ A central
simple algebra is a simple algebra. An algebra $A$ is called \textit{%
alternative} if $x^{2}y=x\left( xy\right) $ and $xy^{2}=\left( xy\right) y,$
for all $x,y\in A,$ \textit{\ flexible} if $x\left( yx\right) =\left(
xy\right) x=xyx,$ for all $x,y\in A$ and \textit{power associative} if the
subalgebra $<x>$ of $A$ generated by any element $x\in A$ is associative. $\ 
$Each alternative algebra is$\ $a\ flexible algebra and a power associative
algebra. In each alternative algebra $A,$ the following identities\newline
1) $a(x(ay))=(axa)y$\newline
2) $((xa)y)a=x(aya)$\newline
3) $(ax)(ya)=a(xy)a$\newline
hold, for all $a,x,y\in A.$ These identities are called the \textit{Moufang
identities.}

A unitary algebra $A\neq K$ such that we have $x^{2}+\alpha _{x}x+\beta
_{x}=0$ for each $x\in A,$ with $\alpha _{x},\beta _{x}\in K,$ is called a 
\textit{quadratic algebra}.

In the following, we briefly present the \textit{Cayley-Dickson process} and
the properties of the algebras obtained. For details about the
Cayley-Dickson process, the reader is referred to $\left[ \text{Sc; 66}%
\right] $ and [Sc; 54].

Let $A$ be a finite dimensional unitary algebra over a field $\ K$ with a 
\textit{scalar} \textit{involution} $\,$%
\begin{equation*}
\,\,\,\overline{\phantom{x}}:A\rightarrow A,a\rightarrow \overline{a},
\end{equation*}%
$\,\,$ i.e. a linear map satisfying the following relations:$\,\,\,\,\,$%
\begin{equation*}
\overline{ab}=\overline{b}\overline{a},\,\overline{\overline{a}}=a,
\end{equation*}%
$\,\,$and 
\begin{equation*}
a+\overline{a},a\overline{a}\in K\cdot 1\ \text{for all }a,b\in A.\text{ }
\end{equation*}%
The element $\,\overline{a}$ is called the \textit{conjugate} of the element 
$a,$ the linear form$\,\,$%
\begin{equation*}
\,\,t:A\rightarrow K\,,\,\,t\left( a\right) =a+\overline{a}
\end{equation*}%
and the quadratic form 
\begin{equation*}
n:A\rightarrow K,\,\,n\left( a\right) =a\overline{a}\ 
\end{equation*}%
are called the \textit{trace} and the \textit{norm \ }of \ the element $a,$
respectively$.$ Hence an algebra $A$ with a scalar involution is quadratic. $%
\,$

Let$\,\,\,\gamma \in K$ \thinspace be a fixed non-zero element. We define
the following algebra multiplication on the vector space 
\begin{equation*}
A\oplus A:\left( a_{1},a_{2}\right) \left( b_{1},b_{2}\right) =\left(
a_{1}b_{1}+\gamma \overline{b_{2}}a_{2},a_{2}\overline{b_{1}}%
+b_{2}a_{1}\right) .
\end{equation*}%
\newline
We obtain an algebra structure over $A\oplus A,$ denoted by $\left( A,\gamma
\right) $ and called the \textit{algebra obtained from }$A$\textit{\ by the
Cayley-Dickson process.} $\,$We have $\dim \left( A,\gamma \right) =2\dim A$.

Let $x\in \left( A,\gamma \right) $, $x=\left( a_{1},a_{2}\right) $. The map 
\begin{equation*}
\,\,\,\overline{\phantom{x}}:\left( A,\gamma \right) \rightarrow \left(
A,\gamma \right) \,,\,\,x\rightarrow \bar{x}\,=\left( \overline{a}_{1},\text{%
-}a_{2}\right) ,
\end{equation*}%
\newline
is a scalar involution of the algebra $\left( A,\gamma \right) $, extending
the involution $\overline{\phantom{x}}\,\,\,$of the algebra $A.$ Let 
\begin{equation*}
\,t\left( x\right) =t(a_{1})
\end{equation*}%
and$\,\,\,$ 
\begin{equation*}
n\left( x\right) =n\left( a_{1}\right) -\gamma n(a_{2})
\end{equation*}%
be $\,\,$the \textit{trace} and the \textit{norm} of the element $x\in $ $%
\left( A,\gamma \right) ,$ respectively.\thinspace $\,$

\thinspace If we take $A=K$ \thinspace and apply this process $t$ times, $%
t\geq 1,\,\,$we obtain an algebra over $K,\,\,$%
\begin{equation}
A_{t}=\left( \frac{\alpha _{1},...,\alpha _{t}}{K}\right) .  \tag{1.1.}
\end{equation}%
By induction in this algebra, the set $\{1,e_{2},...,e_{n}\},n=2^{t},$
generates a basis with the properties:%
\begin{equation}
e_{i}^{2}=\gamma _{i}1,\,\,_{i}\in K,\gamma _{i}\neq 0,\,\,i=2,...,n 
\tag{1.2.}
\end{equation}%
and \ 
\begin{equation}
e_{i}e_{j}=-e_{j}e_{i}=\beta _{ij}e_{k},\,\,\beta _{ij}\in K,\,\,\beta
_{ij}\neq 0,i\neq j,i,j=\,\,2,...n,  \tag{1.3.}
\end{equation}%
$\ \beta _{ij}$ and $e_{k}$ being uniquely determined by $e_{i}$ and $e_{j}.$

From [Sc; 54], Lemma 4, it results that in any algebra $A_{t}$ with the
basis \newline
$\{1,e_{2},...,e_{n}\}$ satisfying relations $\left( 1.2.\right) $ and $%
\left( 1.3.\right) $ we have:

\begin{equation}
e_{i}\left( e_{i}x\right) =\gamma _{i}^{2}=(xe_{i})e_{i},  \tag{1.4.}
\end{equation}%
for all $i\in \{1,2,...,n\}$ and for \ every $x\in A$

It is known that if an algebra $A$ is finite-dimensional, then it is a
division algebra if and only if $A$ does not contain zero divisors. (See
[Sc;66])

Algebras $A_{t}$ of dimension $2^{t}\ $obtained by the Cayley-Dickson
process, described above, are central-simple, flexible and\textit{\ }power
associative for all $t\geq 1$ and, in general, are not division algebras for
all $t\geq 1$.\ But there are fields on which, if we apply the
Cayley-Dickson process, the resulting algebras $A_{t}\ $are division \
algebras for all $t\geq 1.$ (See [Br; 67] and [Fl; 12] ). We remark that the
field $K$ is the center of the algebra $A_{t},\ \ $for $t\geq 2.$(See [Sc;
54])

Let $K$ be a field containing $\omega \,,\ $a primitive $n-$th root of
unity, \ and $A$ be an associative algebra over $K.$ Let $%
S=\{e_{1},...,e_{r}\}$ be a set of elements in $A$ such that the following
condition are fulfilled: $e_{i}e_{j}=\omega e_{j}e_{i}$ for all $i<j$ \ and $%
e_{i}^{n}\in \{1,\omega ,\omega ^{2},...,\omega ^{n-1}\}$. A \textit{%
generalized Clifford algebra} over the field $K,$ denoted by $%
Cl_{r}^{n}\left( K\right) ,$ is defined to be the polynomial algebra $%
K[e_{1},...,e_{r}].$ We remark that the algebra $Cl_{r}^{n}\left( K\right) $
is an associative algebra. For details about generalized Clifford algebra,
the reader is referred to \ [Ki, Ou; 99], [Ko; 10] and [Sm; 91].\medskip

\textbf{Example 1.1.} 1) For $n=2,$ we obtain $Cl_{r}^{2}\left( K\right) $
with $\omega =-1,$ $e_{i}e_{j}=-e_{j}e_{i}$ for all $i<j$ \ and $%
e_{i}^{2}\in \{-1,1\}.$ If\ $\ r=p+q$ and $e_{1}^{2}=...=e_{p}^{2}=1,$ $%
e_{p+1}^{2}=..=e_{q}^{2}=-1,$ then the algebra $Cl_{r}^{2}\left( K\right) $
will be denoted $Cl_{p,q}\left( K\right) .$

2) i) For $p=q=0$ we have $Cl_{0,0}\left( K\right) \simeq K;$

ii) For $p=0,q=1,$ it results that $Cl_{0,1}\left( K\right) $ \ is a
two-dimensional algebra generated by a single vector $e_{1}$ such that $%
e_{1}^{2}=-1$ and therefore $Cl_{0,1}\left( K\right) \simeq K\left(
e_{1}\right) $. For $K=\mathbb{R}$ it follows that $Cl_{0,1}\left( \mathbb{R}%
\right) \simeq \mathbb{C}.$

iii) For $p=0,q=2,$ the algebra $Cl_{0,2}\left( K\right) $ is a
four-dimensional algebra spanned by the set $\{1,e_{1},e_{2},e_{1}e_{2}\}.$
Since $e_{1}^{2}=e_{2}^{2}=(e_{1}e_{2})^{2}=-1$ and $e_{1}e_{2}=-e_{2}e_{1},$
we obtain that this algebra is isomorphic to the division quaternions
algebra $\mathbb{H}$.

iv) For $p=1,q=1$ or $p=2,q=0,$ we obtain the algebra $Cl_{1,1}\left(
K\right) \simeq Cl_{2,0}\left( K\right) $ which is isomorphic with a split
quaternion algebra, called \textit{paraquaternion algebra }or \textit{%
antiquaternion algebra}. (See [Iv, Za; 05])

\begin{equation*}
\end{equation*}

\bigskip

\textbf{2. Main results}%
\begin{equation*}
\end{equation*}

Let $A$ be an algebra obtained by the Cayley-Dickson process with the basis $%
\{e_{0}:=1,e_{1},...,e_{n}\}$ such that,$\ e_{m}e_{r}=-e_{r}e_{m},$ $r\neq
m,e_{m}^{2}=\gamma _{m}\in K,m\in \{1,2,...,n\}.$ For elements $%
a=\sum\limits_{m=0}^{n}a_{m}e_{m},b=\sum\limits_{m=0}^{n}b_{m}e_{m}$ we
define an element in $K$, denoted by \ $T\left( a,b\right) ,$ $T\left(
a,b\right) =\sum\limits_{m=0}^{n}e_{m}^{2}a_{m}b_{m}.\ \ $We denote by $%
\overrightarrow{A}$ the set the elements \thinspace $\{\overrightarrow{a}%
~\mid ~\overrightarrow{a}=\sum\limits_{m=1}^{n}a_{m}e_{m},a_{m}\in K\}.$ It
results that the conjugate of the element $a$ can be written as $\overline{a}%
=a_{0}-~\overrightarrow{a}.$ \ Obviously, $\overrightarrow{\left( ~%
\overrightarrow{a}\right) }=~\overrightarrow{a}$ and ~$\overrightarrow{e_{m}}%
=e_{m}.\bigskip $

\textbf{Lemma 2.1.} \textit{Let }$A$ \textit{be an algebra obtained by the
Cayley-Dickson process. The following equalities are true:}

1) 
\begin{equation*}
T\left( a,b\right) =T\left( b,a\right) ,
\end{equation*}
for all $a,b\in A.$

2) 
\begin{equation*}
T\left( \lambda a,b\right) =\lambda T\left( a,b\right) ,
\end{equation*}
for all $\lambda \in K,~a,b\in A.$

3) 
\begin{equation*}
T\left( a,b+c\right) =T\left( a,b\right) +T\left( a,c\right) ,
\end{equation*}
for all $a,b,c\in A.$

4) 
\begin{equation*}
T\left( a,\overline{a}\right) =a\overline{a}=n\left( a\right) ,
\end{equation*}%
for all $a\in A$

5) \ 
\begin{equation}
\overrightarrow{a}\overrightarrow{b}=2T\left( \overrightarrow{a},%
\overrightarrow{b}\right) -\overrightarrow{b}\overrightarrow{a},  \tag{2.1.}
\end{equation}%
\begin{equation}
ab=ba-2\overrightarrow{b}\overrightarrow{a}+2T\left( \overrightarrow{a},%
\overrightarrow{b}\right) ,  \tag{2.2.}
\end{equation}

\begin{equation}
\overrightarrow{\overrightarrow{a}\overrightarrow{b}}=-T\left( 
\overrightarrow{a},\overrightarrow{b}\right) +\overrightarrow{a}%
\overrightarrow{b}.  \tag{2.3.}
\end{equation}

\begin{equation}
(\overrightarrow{a})^{2}\in K,  \tag{2.4.}
\end{equation}%
for all $a,b\in A.\medskip $

\textbf{Proof. }

Relation from 1), 2), 3), 4) are obvious.

5) For $\overrightarrow{a}=\sum\limits_{m=1}^{n}a_{m}e_{m},\overrightarrow{b}%
=\sum\limits_{m=1}^{n}b_{m}e_{m}$ we obtain

\ 
\begin{equation}
\overrightarrow{a}\overrightarrow{b}\text{=}\sum\limits_{m=1}^{n}a_{m}e_{m}%
\cdot \sum\limits_{m=1}^{n}b_{m}e_{m}\text{=}\sum%
\limits_{m=1}^{n}e_{m}^{2}a_{m}b_{m}\text{+}\alpha \text{=}T\left( 
\overrightarrow{a},\overrightarrow{b}\right) \text{+}\alpha ,\alpha \in 
\overrightarrow{A}.  \tag{2.5.}
\end{equation}%
Computing $\overrightarrow{b}\overrightarrow{a}$, it follows that \ 
\begin{equation}
\overrightarrow{b}\overrightarrow{a}=T\left( \overrightarrow{a},%
\overrightarrow{b}\right) -\alpha ,\alpha \in \overrightarrow{A}.  \tag{2.6.}
\end{equation}%
\newline
If we add relations $\left( 2.5.\right) $ and $\left( 2.6.\right) $, it
results $\overrightarrow{a}\overrightarrow{b}+\overrightarrow{b}%
\overrightarrow{a}=2T\left( \overrightarrow{a},\overrightarrow{b}\right) ,$
therefore relation $\left( 2.1.\right) $ is obtained.

For $a=a_{0}+\overrightarrow{a}$ and $b=b_{0}+\overrightarrow{b},$ we
compute 
\begin{equation*}
ab=\left( a_{0}+\overrightarrow{a}\right) \left( b_{0}+\overrightarrow{b}%
\right) =a_{0}b_{0}+a_{0}\overrightarrow{b}+b_{0}\overrightarrow{a}+%
\overrightarrow{a}\overrightarrow{b}
\end{equation*}%
and%
\begin{equation*}
ba=\left( b_{0}+\overrightarrow{b}\right) \left( a_{0}+\overrightarrow{a}%
\right) =b_{0}a_{0}+b_{0}\overrightarrow{a}+a_{0}\overrightarrow{b}+%
\overrightarrow{b}\overrightarrow{a}.
\end{equation*}

Subtracting the last two relations and using relation $\left( 2.1.\right) $,
we obtain $ab-ba=$ $\overrightarrow{a}\overrightarrow{b}-\overrightarrow{b}%
\overrightarrow{a}=2T\left( \overrightarrow{a},\overrightarrow{b}\right) -2%
\overrightarrow{b}\overrightarrow{a},$ then relation $\left( 2.2.\right) $
is proved.

Relation $\left( 2.3.\right) $ is obvious.

$\bigskip $For $\overrightarrow{a}=\sum\limits_{m=1}^{n}a_{m}e_{m},$ it
results that $(\overrightarrow{a})^{2}=\sum\limits_{m=1}^{n}(a_{m})^{2}\in
K.\Box $

For quaternion algebras, the above result was proved in [Sz; 09].\bigskip

\textbf{Theorem 2.2.} \ \textit{Let} $A$ \textit{be an algebra obtained by
the Cayley-Dickson process such that }$e_{m}^{2}=-1,$\textit{\ for all }$%
m\in \{1,2,...n\}$\textit{. If\ }$n-1\in K-\{0\},$ \textit{then, for \ all} $%
x\in A,$ \textit{we have} 
\begin{equation*}
\overline{x}=\frac{1}{1-n}\underset{m=0}{\overset{n}{\sum }}e_{m}xe_{m}.
\end{equation*}

\textbf{Proof. }Let \ $x=\underset{m=0}{\overset{n}{\sum }}e_{m}x_{m}$. From
Lemma 2.1 and relation $\left( 1.4.\right) ,$ we obtain\newline
$\underset{m=0}{\overset{n}{\sum }}e_{m}xe_{m}=x+\underset{m=1}{\overset{n}{%
\sum }}e_{m}xe_{m}=$ \newline
$=x+\underset{m=1}{\overset{n}{\sum }}e_{m}\left( e_{m}x-2e_{m}%
\overrightarrow{x}+2T\left( e_{m},\overrightarrow{x}\right) \right) =$%
\newline
$=x+\underset{m=1}{\overset{n}{\sum }}e_{m}^{2}x-2\underset{m=1}{\overset{n}{%
\sum }}e_{m}^{2}\overrightarrow{x}+2\underset{m=1}{\overset{n}{\sum }}%
e_{m}^{2}e_{m}x_{m}=$\newline
$=x-nx+2n\overrightarrow{x}-2\underset{m=1}{\overset{n}{\sum }}e_{m}x_{m}=$%
\newline
$=\left( 1-n\right) x-2\left( 1-n\right) \overrightarrow{x}=\left(
1-n\right) \left( x-2\overrightarrow{x}\right) =$\newline
$=\left( 1-n\right) \overline{x}.\Box \medskip $

For the real quaternions, the above relation is well known:

\begin{equation*}
\overline{x}=-\frac{1}{2}\left( x+ixi+jxj+kxk\right) .
\end{equation*}

\textbf{Theorem 2.3. }\textit{\ Let} $A$ \textit{be an algebra obtained by
the Cayley-Dickson process. Then for all} $x,y,z\in A,$ \textit{it results
that} 
\begin{equation}
\left( xy-yx\right) ^{2}z=z\left( xy-yx\right) ^{2}.  \tag{2.7.}
\end{equation}

\textbf{Proof. }

We will compute both members of the equality $\left( xy-yx\right) ^{2}z$=$%
z\left( xy-yx\right) ^{2}.$ Using relation $\left( 2.2.\right) $ from Lemma
1 and since $T\left( \overrightarrow{x},\overrightarrow{y}\right) \in K$, \
we obtain \newline
$\left( -2\overrightarrow{y}\overrightarrow{x}+2T\left( \overrightarrow{x},%
\overrightarrow{y}\right) \right) ^{2}z=z\left( -2\overrightarrow{y}%
\overrightarrow{x}+2T\left( \overrightarrow{x},\overrightarrow{y}\right)
\right) ^{2}\Rightarrow $\newline
$\Rightarrow \left[ 4\left( \overrightarrow{y}\overrightarrow{x}\right)
^{2}+4T^{2}\left( \overrightarrow{x},\overrightarrow{y}\right) -8\left( 
\overrightarrow{y}\overrightarrow{x}\right) T\left( \overrightarrow{x},%
\overrightarrow{y}\right) \right] z=$\newline
$=z\left[ 4\left( \overrightarrow{y}\overrightarrow{x}\right)
^{2}+4T^{2}\left( \overrightarrow{x},\overrightarrow{y}\right) -8\left( 
\overrightarrow{y}\overrightarrow{x}\right) T\left( \overrightarrow{x},%
\overrightarrow{y}\right) \right] \Rightarrow $\newline
$\Rightarrow 4\left( \overrightarrow{y}\overrightarrow{x}\right)
^{2}z+4T^{2}\left( \overrightarrow{x},\overrightarrow{y}\right) z-8T\left( 
\overrightarrow{x},\overrightarrow{y}\right) \left( \overrightarrow{y}%
\overrightarrow{x}\right) z=$\newline
$=4z\left( \overrightarrow{y}\overrightarrow{x}\right) ^{2}+4T^{2}\left( 
\overrightarrow{x},\overrightarrow{y}\right) z-8T\left( \overrightarrow{x},%
\overrightarrow{y}\right) z\left( \overrightarrow{y}\overrightarrow{x}%
\right) .$\newline
Dividing this last relation by $4$ and\ after reducing the terms, it results%
\newline
$\left( \overrightarrow{y}\overrightarrow{x}\right) ^{2}z-2T\left( 
\overrightarrow{x},\overrightarrow{y}\right) \left( \overrightarrow{y}%
\overrightarrow{x}\right) z=z\left( \overrightarrow{y}\overrightarrow{x}%
\right) ^{2}-2T\left( \overrightarrow{x},\overrightarrow{y}\right) z\left( 
\overrightarrow{y}\overrightarrow{x}\right) .$ \newline
We denote 
\begin{eqnarray*}
E &=&\left[ \left( \overrightarrow{y}\overrightarrow{x}\right) ^{2}z-z\left( 
\overrightarrow{y}\overrightarrow{x}\right) ^{2}\right] - \\
&&-\left[ 2T\left( \overrightarrow{x},\overrightarrow{y}\right) \left( 
\overrightarrow{y}\overrightarrow{x}\right) z-2T\left( \overrightarrow{x},%
\overrightarrow{y}\right) z\left( \overrightarrow{y}\overrightarrow{x}%
\right) \right]
\end{eqnarray*}%
and we will prove that $E=0.$\newline
We denote $\ $%
\begin{equation*}
E_{1}=\left( \overrightarrow{y}\overrightarrow{x}\right) ^{2}z-2T\left( 
\overrightarrow{x},\overrightarrow{y}\right) \left( \overrightarrow{y}%
\overrightarrow{x}\right) z
\end{equation*}%
and 
\begin{equation*}
E_{2}=z\left( \overrightarrow{y}\overrightarrow{x}\right) ^{2}-2T\left( 
\overrightarrow{x},\overrightarrow{y}\right) z\left( \overrightarrow{y}%
\overrightarrow{x}\right) .\newline
\end{equation*}%
First, we compute $E_{1}.$ We obtain\newline
$E_{1}=[\left( \overrightarrow{y}\overrightarrow{x}\right) ^{2}-2T\left( 
\overrightarrow{x},\overrightarrow{y}\right) \left( \overrightarrow{y}%
\overrightarrow{x}\right) ]z.$ \newline
From Lemma 2.1., relation $\left( 2.3.\right) ,$ we have $\overrightarrow{y}%
\overrightarrow{x}$= $T\left( \overrightarrow{y},\overrightarrow{x}\right) +$
$\overrightarrow{\overrightarrow{y}\overrightarrow{x}}.$ \newline
Then $\left( \overrightarrow{y}\overrightarrow{x}\right) ^{2}=T^{2}\left( 
\overrightarrow{y},\overrightarrow{x}\right) +\left( \overrightarrow{%
\overrightarrow{y}\overrightarrow{x}}\right) ^{2}+2T\left( \overrightarrow{y}%
,\overrightarrow{x}\right) $ $\overrightarrow{\overrightarrow{y}%
\overrightarrow{x}}.$ \newline
Therefore\newline
$E_{1}$=$[T^{2}\left( \overrightarrow{y},\overrightarrow{x}\right) $+$\left( 
\overrightarrow{\overrightarrow{y}\overrightarrow{x}}\right) ^{2}$+\newline
+$2T\left( \overrightarrow{y},\overrightarrow{x}\right) \overrightarrow{%
\overrightarrow{y}\overrightarrow{x}}$-$2T\left( \overrightarrow{x},%
\overrightarrow{y}\right) \left( \overrightarrow{y}\overrightarrow{x}\right)
]z$=\newline
$=[T^{2}\left( \overrightarrow{y},\overrightarrow{x}\right) $+$\left( 
\overrightarrow{\overrightarrow{y}\overrightarrow{x}}\right) ^{2}$+\newline
+$2T\left( \overrightarrow{y},\overrightarrow{x}\right) (\overrightarrow{%
\overrightarrow{y}\overrightarrow{x}}-\overrightarrow{y}\overrightarrow{x}%
)]z.$\newline
Since $\overrightarrow{\overrightarrow{y}\overrightarrow{x}}-\overrightarrow{%
y}\overrightarrow{x}=-T\left( \overrightarrow{y},\overrightarrow{x}\right) ,$
it results that \newline
$[\left( \overrightarrow{y}\overrightarrow{x}\right) ^{2}-2T\left( 
\overrightarrow{x},\overrightarrow{y}\right) \left( \overrightarrow{y}%
\overrightarrow{x}\right) ]$=\newline
=$[\left( \overrightarrow{\overrightarrow{y}\overrightarrow{x}}\right)
^{2}-T^{2}\left( \overrightarrow{y},\overrightarrow{x}\right) ]$=$\alpha \in
K,$\newline
from Lemma 2.1., relation $\left( 2.4.\right) .$ Hence $E_{1}=\alpha z.$ 
\newline
Now, we compute $E_{2}.$ We obtain\newline
$E_{2}=z[\left( \overrightarrow{y}\overrightarrow{x}\right) ^{2}-2T\left( 
\overrightarrow{x},\overrightarrow{y}\right) \left( \overrightarrow{y}%
\overrightarrow{x}\right) ]=$\newline
$=z\alpha =\alpha z$ since $\alpha \in K.$

It follows that $E=E_{1}-E_{2}=0,$ therefore relation $\left( 2.7.\right) $
is proved. $\Box $\bigskip

\textbf{Remark 2.4.} 1) Identity $\left( 2.7.\right) $ is called the \textit{%
Hall identity}. From the above theorem, we remark that Hall identity is true
for all algebras obtained by the Cayley-Dickson process.

2) Relation $\left( 2.7.\right) $ can be written: $\left[ x,y\right] ^{2}z=z%
\left[ x,y\right] ^{2}$ or $\left[ \left[ x,y\right] ^{2},z\right] =0,$
where $\left[ x,y\right] =xy-yx$ is the commutator of two elements. If $A=%
\mathbb{H},$ then the identity $\left( 2.7.\right) $ is proved by Hall in
[Ha; 43]. If \ $A=\mathbb{H}$ and, for example, $y=i,z=j,$ we have a
quadratic quaternionic equation for which any quaternion is a root:%
\begin{equation*}
xixk\text{+}kxix\text{+}ixixj-jxixi\text{+}x^{2}j-jx^{2}-ix^{2}k-kx^{2}i%
\text{=}0.
\end{equation*}

\textbf{Proposition 2.5.} \textit{Let} $A$ \textit{be an arbitrary algebra
over the field }$K$\textit{\ such that the relation} $\left( 2.7.\right) $%
\textit{\ holds for all} $x,y,z\in A.$ \textit{Then we have }%
\begin{equation}
\left[ \left[ x,y\right] \left[ u,y\right] ,z\right] \text{+}[[x,y][x,v],z]%
\text{+}[[u,y][x,y],z]+[[x,v][x,y],z]\text{=}0,\newline
\tag{2.8.}
\end{equation}%
\begin{equation}
\left[ \lbrack x,v][u,y],z\right] \text{+}\left[ [u,y][x,v],z\right] \text{+}%
\left[ [x,y][u,v],z\right] \text{+}\left[ [u,v][x,y],z\right] \text{=}0,%
\newline
\tag{2.9.}
\end{equation}%
\begin{equation}
\lbrack \lbrack u,y][u,v],z]\text{+}[[x,v][u,v],z]\text{+}[[u,v][u,y],z]%
\text{+}[[u,v][x,v],z]\text{=}0\newline
\tag{2.10.}
\end{equation}%
\textit{for all} $x,y,z,u,v\in A.\medskip $

\textbf{Proof. }$~$We linearize relation $\left( 2.7.\right) .$ Let $%
x,y,z\in A$ be three arbitrary elements such that $\left( xy-yx\right)
^{2}z=z\left( xy-yx\right) ^{2}.$ \newline
For $x+\lambda u,y+\lambda v,z$ we obtain\newline
$[\left( x+\lambda u\right) \left( y+\lambda v\right) -\left( y+\lambda
v\right) \left( x+\lambda u\right) ]^{2}z=$\newline
$=z[\left( x+\lambda u\right) \left( y+\lambda v\right) -\left( y+\lambda
v\right) \left( x+\lambda u\right) ]^{2}.$\newline
It results\newline
$[xy-yx$+$\lambda (uy+xv-yu-vx)$+$\lambda ^{2}\left( uv-vu\right) ]^{2}z=$%
\newline
$=z\left[ xy-yx\text{+}\lambda (uy+xv-yu-vx)\text{+}\lambda ^{2}\left(
uv-vu\right) \right] ^{2}.$ \newline
We obtain \newline
$\left( xy-yx\right) ^{2}z$+$\lambda ^{2}[\left( uy-yu\right) $+$\left(
xv-vx\right) ]^{2}z+$\newline
$+\lambda ^{4}\left( uv-vu\right) ^{2}z+$\newline
$+\lambda \lbrack \left( xy-yx\right) \left( \left( uy-yu\right) \text{+}%
\left( xv-vx\right) \right) ]z+$\newline
$+\lambda \lbrack \left( \left( uy-yu\right) +\left( xv-vx\right) \right)
\left( xy-yx\right) ]z+$\newline
$+\lambda ^{2}[\left( uv-vu\right) \left( xy-yx\right) ]z+$\newline
$+\lambda ^{2}[\left( xy-yx\right) \left( uv-vu\right) ]z+$\newline
$+\lambda ^{3}[[\left( uy-yu\right) $+$\left( xv-vx\right) ]\left(
uv-vu\right) ]z+$\newline
$+\lambda ^{3}[\left( uv-vu\right) [\left( uy-yu\right) $+$\left(
xv-vx\right) ]]z=$\newline
$z\left( xy-yx\right) ^{2}$+$\lambda ^{2}z[\left( uy-yu\right) $+$\left(
xv-vx\right) ]^{2}+$\newline
$+\lambda ^{4}z\left( uv-vu\right) ^{2}+$\newline
$+\lambda z[\left( xy-yx\right) \left( \left( uy-yu\right) \text{+}\left(
xv-vx\right) \right) ]+$\newline
$+\lambda z[\left( \left( uy-yu\right) \text{+}\left( xv-vx\right) \right)
\left( xy-yx\right) ]+$\newline
$+\lambda ^{2}z[\left( uv-vu\right) \left( xy-yx\right) ]+$\newline
$+\lambda ^{2}z[\left( xy-yx\right) \left( uv-vu\right) ]+$\newline
$+\lambda ^{3}z[[\left( uy-yu\right) $+$\left( xv-vx\right) ]\left(
uv-vu\right) ]+$\newline
$+\lambda ^{3}z[\left( uv-vu\right) [\left( uy-yu\right) $+$\left(
xv-vx\right) ]],$ for all $x,y,z,u,v\in A.$ \newline
Since the coefficients of $\lambda $ are equal in both members of the
equality, we obtain:\newline
$[\left( xy-yx\right) \left( \left( uy-yu\right) \text{+}\left( xv-vx\right)
\right) ]z+$\newline
$+[\left( \left( uy-yu\right) \text{+}\left( xv-vx\right) \right) \left(
xy-yx\right) ]z=$\newline
$=z[\left( xy-yx\right) \left( \left( uy-yu\right) \text{+}\left(
xv-vx\right) \right) ]+$\newline
$+z[\left( \left( uy-yu\right) +\left( xv-vx\right) \right) \left(
xy-yx\right) ].$ \newline
We can write this last relation under the form:\newline
$\{\left[ x,y\right] \left[ u,y\right] \}z$+$\{[x,y]\left[ x,v\right] \}z+$%
\newline
$+\{\left[ u,y\right] [x,y]\}z+\{\left[ x,v\right] \left[ x,y\right] \}z=$%
\newline
$=z\{\left[ x,y\right] \left[ u,y\right] \}+z\{[x,y]\left[ x,v\right] \}+$%
\newline
$+z\{\left[ u,y\right] [x,y]\}+z\{\left[ x,v\right] \left[ x,y\right] \}.$ 
\newline
It results\newline
$\left[ \left[ x,y\right] \left[ u,y\right] ,z\right] $+$[[x,y][x,v],z]$+$%
[[u,y][x,y],z]$+$[[x,v][x,y],z]$=$0$\newline
and we obtain relation $\left( 2.8.\right) .$\newline
Since the coefficients of $\lambda ^{2}$ are equal in both members of the
equality, we obtain:\newline
$[\left( uy-yu\right) +\left( xv-vx\right) ]^{2}z+$\newline
$+[\left( uv-vu\right) \left( xy-yx\right) ]z+$\newline
$+[\left( xy-yx\right) \left( uv-vu\right) ]z=$\newline
$=z[\left( uy-yu\right) +\left( xv-vx\right) ]^{2}+$\newline
$+z[\left( uv-vu\right) \left( xy-yx\right) ]+$\newline
$+z[\left( xy-yx\right) \left( uv-vu\right) ].$\newline
It results that\newline
$[\left( uy-yu\right) \left( xv-vx\right) ]z$+$[\left( xv-vx\right) \left(
uy-yu\right) ]z+$\newline
$+[\left( uv-vu\right) \left( xy-yx\right) ]z$+$[\left( xy-yx\right) \left(
uv-vu\right) ]z=$\newline
$z[\left( uy-yu\right) \left( xv-vx\right) ]$+$z[\left( xv-vx\right) \left(
uy-yu\right) ]+$\newline
$+z[\left( uv-vu\right) \left( xy-yx\right) ]$+$z[\left( xy-yx\right) \left(
uv-vu\right) ].$ \newline
We can write this last relation under the form:\newline
$\left[ \left[ x,v\right] \left[ u,y\right] ,z\right] $+$\left[ \left[ u,y%
\right] \left[ x,v\right] ,z\right] $+$\left[ \left[ x,y\right] \left[ u,v%
\right] ,z\right] $+$\left[ \left[ u,v\right] \left[ x,y\right] ,z\right] $=$%
0$ \newline
and we obtain relation $\left( 2.9.\right) .$\newline
Since the coefficients of $\lambda ^{3}$ are equal in both members of the
equality, we obtain:\newline
$[[\left( uy-yu\right) +\left( xv-vx\right) ]\left( uv-vu\right) ]z+$\newline
$+[\left( uv-vu\right) [\left( uy-yu\right) +\left( xv-vx\right) ]]z=$%
\newline
$=z[[\left( uy-yu\right) +\left( xv-vx\right) ]\left( uv-vu\right) ]+$%
\newline
$+z[\left( uv-vu\right) [\left( uy-yu\right) +\left( xv-vx\right) ]].$ 
\newline
We can write this last relation under the form:\newline
$[[u,y][u,v],z]$+$[[x,v][u,v],z]$+$[[u,v][u,y],z]$+$[[u,v][x,v],z]$=$0$ 
\newline
and we obtain relation $\left( 2.10.\right) .$ $\Box \bigskip $

\textbf{Remark 2.6.} 1) In [Ti; 99] and [Fl; 01] \ some equations over
division quaternion\ algebra and octonion algebra are solved: in [Fl; 01]
for general case, when $K$ is a commutative field with $charK\neq 2$ and $%
\gamma _{m}$ are arbitrary and in [Ti; 99] for $K=\mathbb{R},$ $\gamma
_{m}=-1,$ with $m\in \{1,2\}$ \ for quaternions and $m\in \{1,2,3\}$ for
octonions. Let $A$ be such an algebra. For example, equation 
\begin{equation}
ax=xb,a,b,x\in A,  \tag{2.11.}
\end{equation}%
for $a\neq \overline{b}$ has general solution under the form $x=%
\overrightarrow{a}p+p\overrightarrow{b},$ for arbitrary $p\in A.$

2) In [Fl, \c{S}t; 09], authors studied equation $x^{2}a=bx^{2}+c,a,b,c\in
A, $ where $A$ is a generalized quaternion division \ algebra or an
generalized octonion division algebra. If $A$ is an arbitrary algebra
obtained by the Cayley-Dickson process and $\ a,b,c\in A$ with $a=b$ and $%
c=0,$ then, from Theorem 2.3., it results that this equation has infinity of
solutions of the form $x=vw-wv,$ where $v,w\in A.\bigskip $

\textbf{Proposition 2.7.} \textit{Let} $A$ \textit{be a quaternion algebra
or an octonion algebra. Then for all} $x,y\in A,$ \textit{there are the
elements} $z,w$ \textit{such that} $(xy-yx)^{2}=\overrightarrow{z}w+w%
\overrightarrow{z}.\medskip $

\textbf{Proof.} Let $z$ be an arbitrary element in $A-K.$ From Theorem 2.3.,
we have that $\left( xy-yx\right) ^{2}z=z\left( xy-yx\right) ^{2},$ for all $%
x,y,z\in A.$ Since $z\neq \overline{z}$ and \ $\left( xy-yx\right) ^{2}$ is
a solution for the equation $\left( 2.11.\right) ,$ from Remark 2.6., it
results that there is an element $w\in A$ such that $(xy-yx)^{2}=%
\overrightarrow{z}w+w\overrightarrow{z}.\Box \bigskip $

\textbf{Proposition 2.8.} \textit{Let }$A$ \textit{be a finite dimensional
unitary algebra with a scalar involution} 
\begin{equation*}
\,\,\,\overline{\phantom{x}}:A\rightarrow A,a\rightarrow \overline{a},
\end{equation*}%
\textit{such that for all} $x,y\in A,$ \textit{the following equality holds:}

\begin{equation}
(x\overline{y}+y\overline{x})^{2}=4\left( x\overline{x}\right) \left( y%
\overline{y}\right) .  \tag{2.12.}
\end{equation}

\textit{Then the algebra} $A$ \textit{has dimension }$1.\medskip $

\textbf{Proof.} \ We remark that $x\overline{y}+y\overline{x}=x\overline{y}+%
\overline{x\overline{y}}\in K.$ First, we prove that \ $[x\overline{y}+y%
\overline{x}]^{2}=4\left( x\overline{x}\right) \left( y\overline{y}\right)
,\forall x,y\in A,$ if and only if $x=ry,$ $r\in K.$ If $x=ry,$ then
relation $\left( 2.12.\right) $ is proved.\ Conversely, assuming that
relation $\left( 2.12.\right) $ is true and supposing that there is not an
element $r\in K$ such that $x=ry,$ then for each two non zero elements $%
a,b\in K,$ we have $ax+by\neq 0.$ Indeed, if $ax+by=0,$ it results $x=-\frac{%
b}{a}y,$ false. We obtain that 
\begin{equation}
\left( ax+by\right) \overline{\left( ax+by\right) }\neq 0.  \tag{2.13.}
\end{equation}%
\newline
Computing relation $\left( 2.13.\right) ,$ it follows 
\begin{equation}
a^{2}\left( x\overline{x}\right) +abx\overline{y}+bay\overline{x}+b^{2}y%
\overline{y}\neq 0.  \tag{2.14.}
\end{equation}%
If we \ put $a=y\overline{y}$ in relation $\left( 2.14.\right) $ and then
simplify by $a,$ it results 
\begin{equation}
\left( y\overline{y}\right) \left( x\overline{x}\right) +bx\overline{y}+by%
\overline{x}+b^{2}\neq 0.  \tag{2.15.}
\end{equation}%
Let $b=-\frac{1}{2}\left( x\overline{y}+y\overline{x}\right) \in K,b\neq 0.$
If we replace this value in relation $\left( 2.15.\right) ,$ we obtain $%
4\left( x\overline{x}\right) \left( y\overline{y}\right) -(x\overline{y}+y%
\overline{x})^{2}\neq 0,$ which it is false. Therefore, there is an element $%
r\in K$ such that $x=ry.$

Assuming that the algebra $A$ has dimension greater or equal with $2,$ it
results that there are two linearly independent vectors, $v$ and $w,$
respectively$.$ Since relation $\left( 2.12.\right) $ is satisfies for $v$
and $w$, we obtain that there is an element $s\in K$ such that $v=sw,$ which
it is false. Hence $\dim A=1.\Box \bigskip $

\textbf{Proposition 2.9.} \textit{Let} $A$ \textit{be an alternative
algebra\ over the field} $K$ \textit{whose center is} $K.$ \textit{If} $%
\left( xy-yx\right) ^{2}z=z\left( xy-yx\right) ^{2}$ \textit{\ for all} $%
x,y,z\in A,$ \textit{then} $A$ \textit{is a quadratic algebra}.\medskip

\textbf{Proof.} Let $x,y\in A$ such that $xy\neq yx.$ If we denote $z=xy-yx,$
it follows that $z^{2}$ commutes with all elements from $A,$ then $z^{2}$ is
in the center of $A.$ \ We obtain $z^{2}=\alpha \in K^{\ast }.$ For $%
t=x^{2}y-yx^{2}$ it results that $t^{2}=(x^{2}y-yx^{2})^{2}\in K$ and $%
t=\left( xy-yx\right) x+x\left( xy-yx\right) =zx+xz.$ We have $zt=z\left(
zx+xz\right) =z^{2}x+zxz=\alpha x+zxz$ and $tz=\left( zx+xz\right)
z=zxz+xz^{2}=\alpha x+zxz.$ Therefore $tz=zt.$ For $z+t=\left(
x^{2}+x\right) y-y\left( x^{2}+x\right) $ $\ $\ we have that $\left(
z+t\right) ^{2}=\beta \in K,$ then $z^{2}+t^{2}+2tz=\beta ,$ hence $%
tz=\gamma \in K.$ Since $\ zx=x(yx)-\left( yx\right) x,$ it follows that $%
(zx)^{2}=\delta \in K.$ If we multiply the relation $\left( zx\right) \left(
zx\right) =\delta $ with $z$ in the left side, we obtain $z\left( \left(
zx\right) \left( zx\right) \right) =\delta z.$ Using alternativity and then
flexibility, it results $\left( z^{2}x\right) \left( zx\right) =\delta z,$
therefore $\ \ \alpha \left( xzx\right) =\delta z,$ hence $xzx=\theta z,$
where $\theta =\alpha ^{-1}\delta .$ It follows that $z\left( xzx\right)
=\theta z^{2}=\theta \alpha \in K.$ Since $z\left( xzx\right) =\left(
zxz\right) x,$ from Moufang identities, we have that $\left( zxz\right)
x=\theta \alpha \in K.$ It results that $\gamma x=\left( tz\right) x=\left(
\alpha x+zxz\right) x=\alpha x^{2}+\left( zxz\right) x=\alpha x^{2}+\theta
\alpha ,$ hence $x^{2}=ax+b,$ where $a=\alpha ^{-1}\gamma ,b=-\theta .$ We
obtain that $A$ is a quadratic algebra.$\Box $

When $A$ is a division associative algebra, this proposition was proved \ by
Hall in \ [Ha; 43].\bigskip\ 

\textbf{Theorem 2.10. }\textit{Let} $\ A$ \textit{be a alternative algebra
such that the center of }$A$ \textit{is }$K$ \textit{and } $\left(
xy-yx\right) ^{2}z=z\left( xy-yx\right) ^{2},$ \textit{\ for all} $x,y,z\in
A.$

1) \textit{If} $A$ \textit{is a division algebra}, \textit{then} $A=K$ 
\textit{or} $A=A_{t},t\in \{1,2,3\}$, \textit{where} $A_{t}$ \textit{is a
division algebra obtained by the Cayley-Dickson process.}

2) \textit{If} $A$ \textit{is not a division algebra } \textit{and} \textit{%
there are two elements} $y,z\in A$ \textit{such that} $y^{2},z^{2}\in K,\
yz=-zy,$ \textit{then }$A$ \textit{is} \textit{a generalized split
quaternion algebra} \textit{\ or }$A$ \textit{is} \textit{a generalized
split octonion algebra.\medskip }

\textbf{Proof.} \ 1) From Proposition 2.9., it results that $A$ is a
quadratic algebra, therefore, from [Al; 49], Theorem 1, we have $\dim A\in
\{ $ $1,2,4,8\}.$ If $\dim A=1,$ then $A=K.$ If $\dim A=2,$ since the \
center is $K,$ then we can find an element $x\in A-K$ such that $x^{2}\in K.$
It results that the set $\{1,x\}$ is a basis in $A,$ therefore $A=K\left(
x\right) $ is a quadratic field extension of the field $K.$ If $\dim A=4,$
from [Al; 39], p. 145, we have that there are two elements $x,y\in A$ such
that $x^{2}=x+a$ ~with $\ 4a+1\neq 0,$ $xy=y\left( 1-x\right)
,y^{2}=b,a,b\in K.$ Denoting $z=x-\frac{1}{2},$ we obtain that $z^{2}=\left(
x-\frac{1}{2}\right) ^{2}=a-\frac{1}{4}\in K.$ and $\ zy=-yz.$ Since $%
zy=\left( x-\frac{1}{2}\right) y=xy-\frac{y}{2}=y-yx-\frac{y}{2}=\frac{y}{2}%
-yx$ and $yz=y\left( x-\frac{1}{2}\right) =yx-\frac{1}{2},$ we have $yz=-zy$
\ then $\left( yz\right) ^{2}\in K.$ It follows that in \ the algebra $A$ we
can find the elements $\ y.z$ such that $y^{2},z^{2},\left( yz\right)
^{2}\in K$ and $yz=-zy.$ Therefore, from [Al; 49], Lemma 4, it results that $%
A$ is a generalized division quaternion algebra. If $A$ has dimension $8$,
denoting with $Q$ the algebra $Q=K+yK+zK+yzK,$ from [Al; 49], Lemma 3, Lemma
4 and Lemma 5, it results that there are the elements $w\in A-Q,$ such that $%
w^{2},\left( yw\right) ^{2},\left( zw\right) ^{2}\in K,yw=-wy,zw=-wz,yz=-zy.$
It follows that $A=K+yK+zK+yzK+wK+wyK+wzK+w(yz)K\ $ is a generalized
division octonion algebra.

2) From the above, it results that $A=Q=K+yK+zK+yzK$ is a generalized
quaternion algebra, which is split from hypothesis. or there is an element $%
w\in A-Q$ such that $A=K+yK+zK+yzK+wK+wyK+wzK+w(yz)K.\ $ In the last case, $%
A $ is a generalized split octonion algebra.$\Box \bigskip $

\textbf{Corollary 2.11. } \textit{Let} $\ A$ \textit{be a non-division
associative algebra such that the center of }$A$ \textit{is }$K.$ \textit{If
\ in algebra }$A$ \textit{we have} $\ \left( xy-yx\right) ^{2}z=z\left(
xy-yx\right) ^{2}$ \textit{\ for all} $x,y,z\in A$ \textit{and} \textit{%
there are two elements} $v,w$ \textit{such that} $v^{2},w^{2}\in K,$ $%
vw=-wv, $ \textit{then }$A$ \textit{is} \textit{a generalized split
quaternion algebra.}$\Box \bigskip $

\textbf{Example 2.12. }

1) Using notations given \ in Preliminaries, if in Theorem 2.10., we have $%
t=1$ and $\alpha _{1}=-1,$ it results that $A=Cl_{0,1}\left( K\right) $ is a
quadratic field extension of the field $K.$ If $t=2$ and $\alpha _{1}=\alpha
_{2}=-1,$ we have that $A=Cl_{0,2}\left( K\right) $ is a quadratic division
quaternion algebra$.$

2) If we have $v^{2},w^{2}\in \{-1,1\}$ in Corollary 2.11$,$ then $%
A=Cl_{1,1}\left( K\right) \simeq Cl_{2,0}\left( K\right) .$%
\begin{equation*}
\end{equation*}

\textbf{Conclusions.} In this paper we proved that the Hall identity is true
in all algebras obtained by the Cayley-Dickson process and that the converse
is also true, in some particular conditions. As we can see in Remark 2.6.,
some identities in algebras obtained by the Cayley-Dickson process can be
used to find solutions for some equations in these algebras or to solve
them. This idea can constitute the starting point for a further
research.\bigskip\ 

\begin{equation*}
\end{equation*}%
\begin{equation*}
\end{equation*}

\textbf{References}%
\begin{equation*}
\end{equation*}

[Al; 39] Albert, A. A., \textit{Structure of algebras}, Amer. Math. Soc.
Colloquium Publications, vol. \textbf{24}, 1939

[Al; 49] Albert, A. A., \textit{Absolute-valued algebraic algebras}, Bull.
Amer. Math. Soc., \textbf{55}(1949), 763-768.

[Br; 67] Brown, R. B., \textit{On generalized Cayley-Dickson algebras},
Pacific J. of Math.,\textbf{\ 20(3)}(1967), 415-422.

[Fl; 12] Flaut, C., \textit{Levels and sublevels of \ algebras obtained by
the Cayley-Dickson process}, 2011, submitted.

[Fl; 01] Flaut, C., \textit{Some equations in algebras obtained by the
Cayley-Dickson process}, An. St. Univ. Ovidius Constanta, 9(2)(2001), 45-68.

[Fl, \c{S}t; 09] Flaut, C., \ \c{S}tef\u{a}nescu, M., \textit{Some equations
over generalized quaternion and octonion division algebras}, Bull. Math.
Soc. Sci. Math. Roumanie, \textbf{52(4)}(100), 2009, 427--439.

[Ha; 43] Hall, M., \textit{Projective planes}, Trans. Amer. Math. Soc. vol. 
\textbf{54}(1943), 229-277.

[Iv, Za; 05] Ivanov, S. Zamkovoy, S., \textit{Parahermitian and
paraquaternionic manifolds} , Differential Geometry and its Applications 
\textbf{23}(2005), 205--234

[Le; 06 ] Lewis, D. W., \textit{Quaternion Algebras and the Algebraic Legacy
of Hamilton's Quaternions}, Irish Math. Soc. Bulletin \textbf{57}(2006),
41--64.

[Ki, Ou; 99] \ El Kinani, E. H., Ouarab, A., \textit{The Embedding of} $%
U_{q}(sl\left( 2\right) )$ \textit{and Sine Algebras in Generalized Clifford
Algebras}, Adv. Appl. Clifford Algebr., \textbf{9(1)}(1999), 103-108.

[Ko; 10] Ko\c{c}, C., \textit{C-lattices and decompositions of generalized
Clifford algebras,} Adv. Appl. Clifford Algebr., \textbf{20(2)}(2010),
313-320.

[Sc; 66] Schafer, R. D., \textit{An Introduction to Nonassociative Algebras,}
Academic Press, New-York, 1966.

[Sc; 54] Schafer, R. D., \textit{On the algebras formed by the
Cayley-Dickson process,} Amer. J. Math., \textbf{76}(1954), 435-446.

[Smi; 50] Smiley, M. F., \textit{A remark on a theorem of Marshall Hall,}
Proceedings of the American Mathematical Society, \textbf{1}(1950), 342-343.

[Sm; 91] Smith T. L., \ \textit{Decomposition of Generalized Clifford
Algebras}, Quart. J. Math. Oxford, 42 (1991), pp. 105-112.

[Sz; 09] Szpakowski, V. S., \textit{Solution of general quadratic
quaternionic equations}, Bull. Soc. Sci. Lettres \L \'{o}d\'{z} 59, Ser.
Rech. D\'{e}form. \textbf{58}(2009), 45 -- 58.

[Ti; 99] Tian, Y., \textit{Similarity and cosimilarity of elements in the
real Cayley-Dickson algebras}, Adv. Appl. Clifford Algebras, 9(1)(1999),
61-76.

\bigskip \bigskip

\bigskip Cristina FLAUT

{\small Faculty of Mathematics and Computer Science,}

{\small Ovidius University,}

{\small Bd. Mamaia 124, 900527, CONSTANTA,}

{\small ROMANIA}

{\small http://cristinaflaut.wikispaces.com/}

{\small http://www.univ-ovidius.ro/math/}

{\small e-mail:}

{\small cflaut@univ-ovidius.ro}

{\small cristina\_flaut@yahoo.com}%
\begin{equation*}
\end{equation*}

Vitalii \ SHPAKIVSKYI

{\small Department of Complex Analysis and Potential Theory}

{\small \ Institute of Mathematics of the National Academy of Sciences of
Ukraine,}

{\small \ 3, Tereshchenkivs'ka st.}

{\small \ 01601 Kiev-4}

{\small \ UKRAINE}

{\small \ http://www.imath.kiev.ua/\symbol{126}complex/}

{\small \ e-mail: shpakivskyi@mail.ru}

\end{document}